\documentclass[twocolumn]{article}
\usepackage{imr}
\usepackage{algorithm}
\usepackage[noend]{algpseudocode}
\usepackage{algorithm} 
\usepackage{algpseudocode}
\usepackage{soul,color}
\algnewcommand\algorithmicforeach{\textbf{for each}}
\algdef{S}[FOR]{ForEach}[1]{\algorithmicforeach\ #1\ \algorithmicdo}
\usepackage{graphicx}
\usepackage{amsfonts}
\usepackage{amssymb}
\usepackage{amsmath}
\usepackage{amsthm}
\usepackage{amsmath, amsthm, amssymb,bm}
\usepackage{comment}

\begin{document}

\title{\uppercase{Towards Tangled Finite Element Analysis over  Partially Inverted Hexahedral Elements}}
\author{Bhagyashree Prabhune$^1$ \and Krishnan Suresh$^2$}
\date{
	University of Wisconsin-Madison, Madison, WI, U.S.A.\\
$^1$bprabhune@wisc.edu \;\; $^2$ksuresh@wisc.edu
}

\abstract{If a finite element mesh contains concave elements, it is said to tangled. Tangled meshes can occur during mesh generation, mesh optimization, and large deformation simulations, and will lead to erroneous results during finite element analysis. Recently, the authors introduced the \emph {tangled finite element method} (TFEM) to accurately handle tangled 2D  concave quadrilateral elements. In this paper, TFEM is extended to 3D hexahedral elements. We demonstrate that TFEM leads to accurate results over tangled hexahedral elements, requiring minimal changes to existing FEM framework.}


\maketitle
\thispagestyle{empty}
\pagestyle{empty}

\section{Introduction}
Hexahedral elements often are preferred over tetrahedral elements for finite element analysis \cite{shepherd2008hexahedral}. However, automatic generation of high quality hexahedral meshes for complex geometries has been evasive.  A critical requirement for a valid finite element mesh is that all elements should be convex \cite{shepherd2008hexahedral}. In other words, the determinant of Jacobian associated with parametric mapping should be positive throughout the element \cite{zienkiewicz2005finite, frey2007mesh, lo2014finite}. If the Jacobian is partly negative, the element is said to be partially-inverted (i.e., concave) and the mesh is said to be tangled. It is well known that finite element analysis over a tangled mesh will lead to erroneous results \cite{zienkiewicz2005finite, prabhune2022tangled}. 

Tangling can occur during mesh generation \cite{gregson2011all, huang2011boundary, nieser2011cubecover}, mesh optimization \cite{xu2018hexahedral}, large deformation simulations \cite{vachal2004untangling}, shape optimization \cite{knupp2003method, staten2011comparison}. To address concave elements, untangling techniques have been proposed \cite{knupp2001hexahedral, knupp2001algebraic, knupp2003method, ruiz2015simultaneous, livesu2015practical, xu2018hexahedral, garanzha2021local}. However, untangling is not always possible/reliable. Other strategies to handle concave elements include advanced FEM techniques like SFEM \cite{liu2007smoothed}, PolyFEM \cite{manzini2014new}, VEM \cite{beirao2014hitchhiker}. However they entail major changes in the existing FEM framework and do not simplify to the standard FEM when the mesh is not tangled.

Recently, the authors proposed the \emph {tangled finite element method} (TFEM) to accurately handle tangled 2D  concave quadrilateral elements \cite{prabhune2022tangled}. The basic concept in TFEM is to introduce certain correction terms to the stiffness matrix when a mesh is tangled, leading to accurate solutions over even severely tangled meshes. In this paper, TFEM is extended to 3D hexahedral meshes. 

The remainder of the paper is organized as follows. Section~\ref{Sec2} describes  the proposed TFEM for 3D hexahedral  meshes. This is followed by  numerical examples in Section~\ref{Sec3}. Conclusion and future work are discussed in Section~\ref{Sec4}.

\section{Tangled Finite Element Method (TFEM)} 
\label{Sec2}

\begin{figure}
	\begin{center}
		\includegraphics[width=0.21\textwidth]{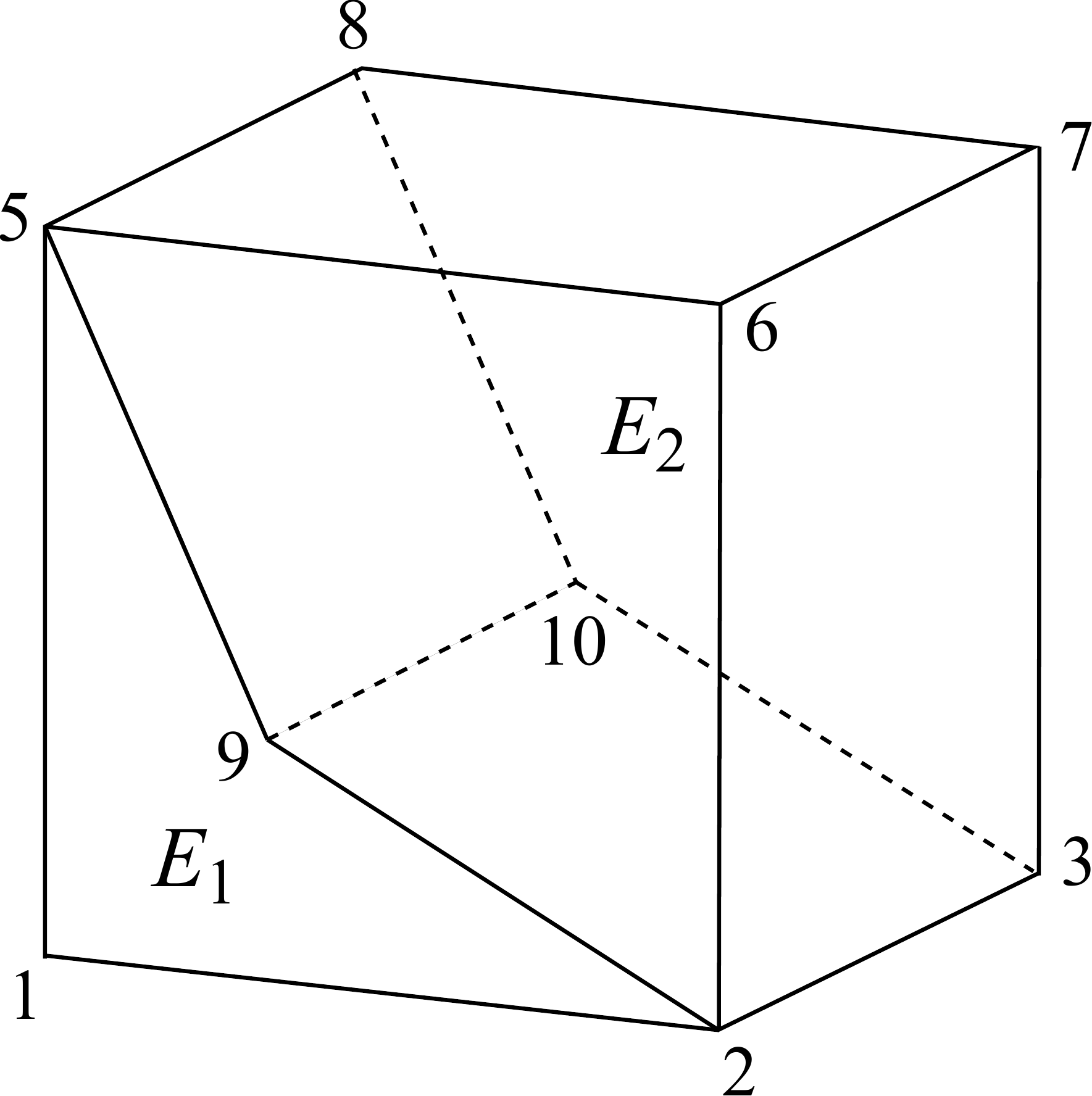}
	\end{center} 	
	\caption{Two-element mesh with a concave element.}
	\label{Fig_MeshHex2elem}
\end{figure}

In this section, we discuss the concept behind TFEM. Consider a two-element hexahedral mesh with one concave element as shown in Fig.~\ref{Fig_MeshHex2elem}. Before discussing the TFEM method, we elaborate on why the concave element leads to erroneous results. 

\begin{figure}
	\begin{center}
		\includegraphics[width=0.48\textwidth]{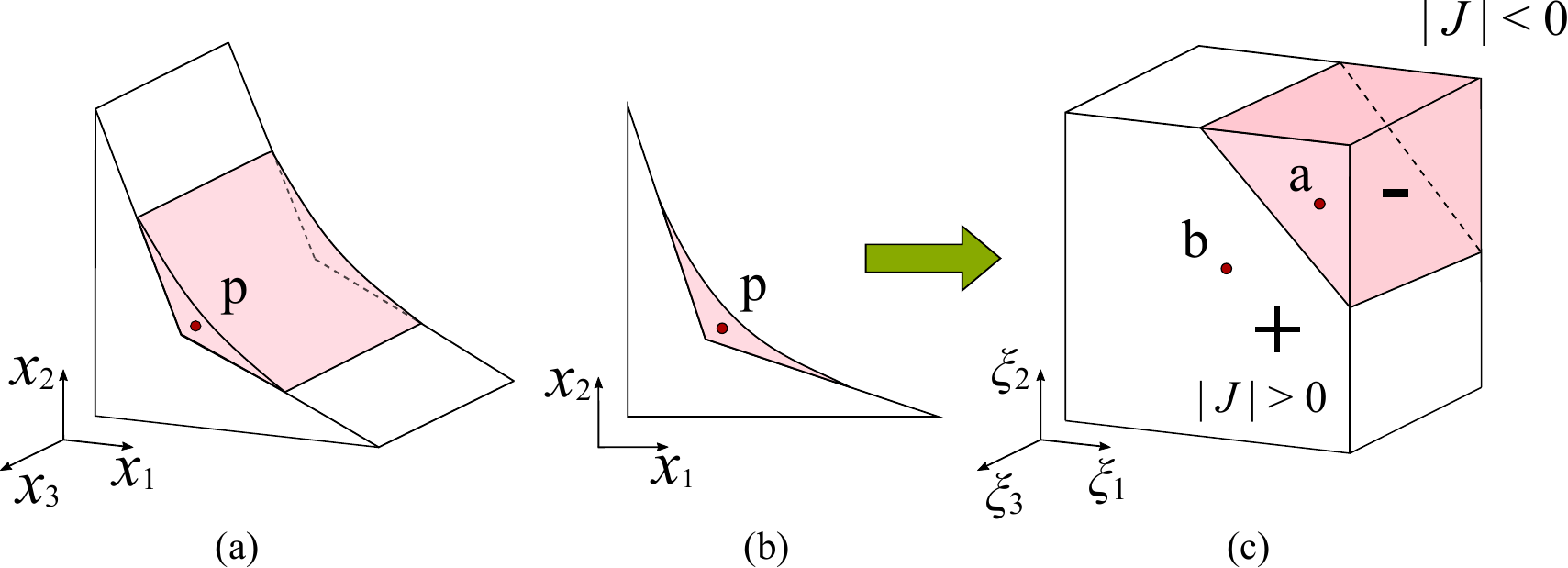}
	\end{center}	
	\caption{(a) Concave hex element.(b) Its projection on $x_1-x_2$ plane. (c) Parametric space.}
	\label{Fig_3D_concaveElem}
\end{figure}

The concave element is shown separately in Fig.~\ref{Fig_3D_concaveElem}a along with its projected view on the $x_1-x_2$ plane in Fig.~\ref{Fig_3D_concaveElem}b.  Further,  Fig.~\ref{Fig_3D_concaveElem}c illustrates the parametric space $(\xi_1,\xi_2, \xi_3)$ associated with this element. Due to the concavity, one can show that the determinant of Jacobian $|\boldsymbol{J}|$ associated with the mapping from the parametric space to the physical space changes sign, i.e., there exists a surface where $|\boldsymbol{J}|$ vanishes, dividing the parametric space into a positive $|\boldsymbol{J}|$ region and a negative $|\boldsymbol{J}|$ region. The negative $|\boldsymbol{J}|$ region in parametric space is highlighted. 
The corresponding set of points in the physical space is also highlighted; this region will be referred as a \emph{tangled} or \emph{folded} region.

Due to the fold, parametric points such as $a$ and $b$ map to the same physical point $p$ in the folded  region. In TFEM, the two parametric regions with positive and negative Jacobians are treated separately and are termed respectively as \textit{positive component} and \textit{negative component}. 

 \begin{figure}
	\begin{center}
		\includegraphics[width=0.48\textwidth]{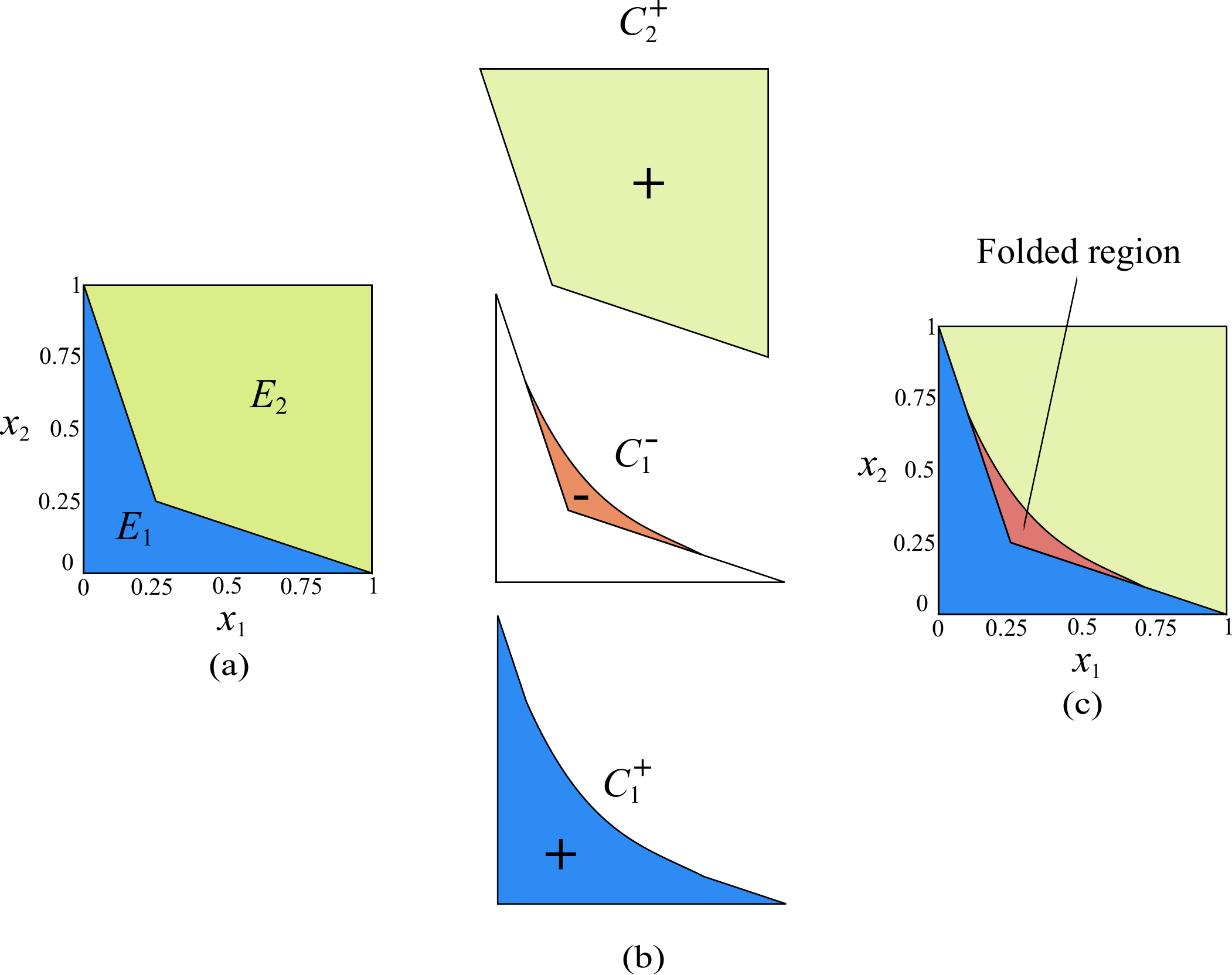}
	\end{center}
	
	\caption{(a) Projection of two-element hex mesh on $x_1-x_2$ plane (b) Components (c) Final physical space is self-overlapping }
	\label{Fig_3D_componentsSeparate}
\end{figure}
 Observe that the concave element $E_1$ can be expressed as the difference between a positive component and a negative component, i.e., $E_1 = C_1^+ - C_1^-$ (see Fig.~\ref{Fig_3D_componentsSeparate}b).  The overlapping region (fold) in Fig.~\ref{Fig_3D_componentsSeparate}c is simply $C_1^-$.

Now to construct the stiffness matrix,  we define two shape functions  $\boldsymbol{N}_1^+$ and $\boldsymbol{N}_1^-$ corresponding to  $C_1^+$ and $C_1^-$ respectively. For example, for the point  $p = (x_1 =5/16,\; x_2 =5/16, \; x_3 = 1)$ in Fig.~\ref{Fig_3D_concaveElem}a, $\boldsymbol{N}_1^+(p)$, are the shape functions of element $E_1$ evaluated at parametric point $b(\xi_1=0, \;\xi_2 = 0,\; \xi_3 = 1)$ whereas $\boldsymbol{N}_1^-(p)$ are the shape functions of element $E_1$ evaluated at point $a (\xi_1 =2/3,\; \xi_2 = 2/3,\; \xi_3 = 1)$. For the convex element, only $\boldsymbol{N}_2^+$ exists, while $\boldsymbol{N}_2^-$ is defined to be zero.

 The critical idea in TFEM is to define the field at a point $p$ within a fold as the oriented sum of the contributions from all \emph{components} the point belongs to. 
 For the above example, since $p$ belongs to three different components $C_1^+, C_1^- $ and $C_2^+$, we have:
\begin{equation}
	u^h(p) = \boldsymbol{N}_1^+(p)\boldsymbol{{\hat u}}_1 - \boldsymbol{N}_1^-(p)  \boldsymbol{{\hat u}}_1 + \boldsymbol{N}_2^+(p) \boldsymbol{{\hat u}}_2
	\label{Eq_fieldAtq}
\end{equation}	
 where $\boldsymbol{{\hat u}}_1$ are the four nodal values of $E_1$,  and $\boldsymbol{{\hat u}}_2$ are the four nodal values of $E_2$. 
  
In addition, TFEM requires an \emph {equality} constraint that ensures uniqueness of the field, and optimal convergence  \cite{prabhune2022tangled}. For this example, equality condition can be stated as:
\begin{equation}
	\boldsymbol{N}_1^+(p)\boldsymbol{{\hat u}}_1  - \boldsymbol{N}_1^-(p)  \boldsymbol{{\hat u}}_1 = 0.
	\label{Eq_equalityCondition}
\end{equation}

If we substitute Eq.~\ref{Eq_equalityCondition} in Eq.~\ref{Eq_fieldAtq}, the field in the tangled region reduces to:
\begin{equation}
u^h(p)  = \boldsymbol{N}_2^+(p) \boldsymbol{{\hat u}}_2
\end{equation}
One can show that, when this is substituted in the standard Galerkin formulation results in the stiffness matrix:
\begin{equation}
	\boldsymbol{K} = 
	\boldsymbol{K}_{convex}^0 +  \boldsymbol{K}_{concave} 
	\label{Eq_FinalImplicit}	
\end{equation}	
where, $\boldsymbol{K}_{convex}^0$ is the stiffness matrix corresponding to the convex element and can be computed by standard FEM procedure. $\boldsymbol{K}_{concave}$ is the stiffness matrix corresponding to the region enclosed within the physical boundary of the concave element (excluding the tangled region) and can be computed by tetrahedralizing this region.  

In summary, TFEM results in the system of equations:
\begin{equation}
	\begin{bmatrix}
		\boldsymbol{K} & \boldsymbol{\tilde{C}}^\top\\
		\boldsymbol{\tilde{C}}& \boldsymbol{0}
	\end{bmatrix} 	
	\begin{Bmatrix}
		\boldsymbol{\hat{u}}\\
		\boldsymbol{\mu}
	\end{Bmatrix} 	= 
	\begin{Bmatrix}
		\boldsymbol{f}^\theta\\
		\boldsymbol{0}
	\end{Bmatrix}.
	\label{Eq_GeneralizedimplicitTFEM}
\end{equation}
where $\boldsymbol{\mu}$ and  $\boldsymbol{\tilde{C}}$ arises from the equality condition. 

\section{Numerical Experiments} 
\label{Sec3}
In this section, we demonstrate the proposed method using numerical experiments. We consider Laplace problem over a unit cube domain $\Omega \in \left( 0,1\right) ^3 $. The field $u$ satisfies the Laplace equation with $u = 0.579x + 0.246y + 0.482z - 0.374$ as the exact solution. The corresponding Dirichlet boundary condition  is applied on left surface while  Neumann condition is applied on the remaining surfaces. The problem is solved over various implicitly tangled meshes. The accuracy is measured via the relative error = $|u_{exact} - u_{computed}|/u_{exact}$ over all nodes. To compute $\boldsymbol{K}_{concave}$, tetrahedralization of the concave element was carried out using Tetgen software. Since Tetgen requires the boundary representation of the region, the boundary of the region is approximately represented as triangles. The accuracy of TFEM depends on the accuracy of triangulations to represent the actual boundary. 

\subsection{Two-element mesh}
\begin{figure}
	
	\begin{center}
		\includegraphics[width=0.37\textwidth]{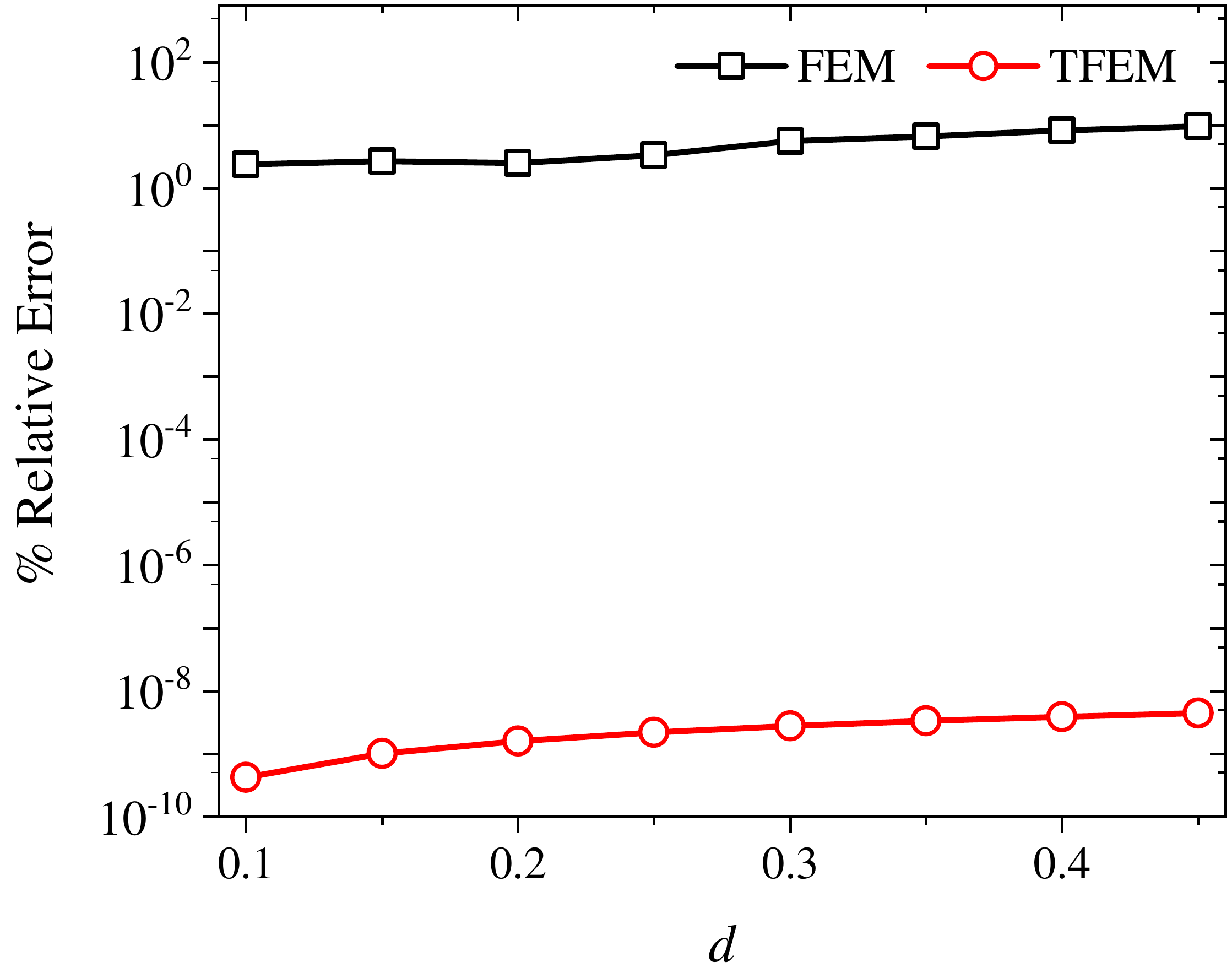}
	\end{center}
	
	\caption{TFEM and FEM errors for two-element mesh.}
	\label{Fig_3D_2Elem_accuracy}
\end{figure}

In Fig.~\ref{Fig_MeshHex2elem}, the domain is divided into two hexahedral elements, one of which is concave. The concave element has two re-entrant vertices: 9 and 10. To vary the extent of tangling, the position of vertex 9 is defined as $ \left( 0.5-d, 0.48-d, 0\right) $  where $ d$ ranges from 0.1 to 0.45. The position of vertex 10 is set to  $ \left( 0.45, 0.4, 1 \right) $. The position of re-entrant vertices is chosen such that the mesh is asymmetric to avoid being simply an extruded 2D mesh. FEM error is around $10^9$ times higher than the TFEM error as seen in Fig.~\ref{Fig_3D_2Elem_accuracy}. The factor preventing the TFEM solution to reach machine precision is the error introduced due to approximating the boundary of complex concave region by finite number of triangles. As the surface representation becomes closer to the actual surface by increasing the number of triangles, the TFEM error decreases. However, it can become computationally infeasible as the number of triangles increases.


\subsection{Eight-element mesh}
\begin{figure}
	
	\begin{center}
		\includegraphics[width=0.35\textwidth]{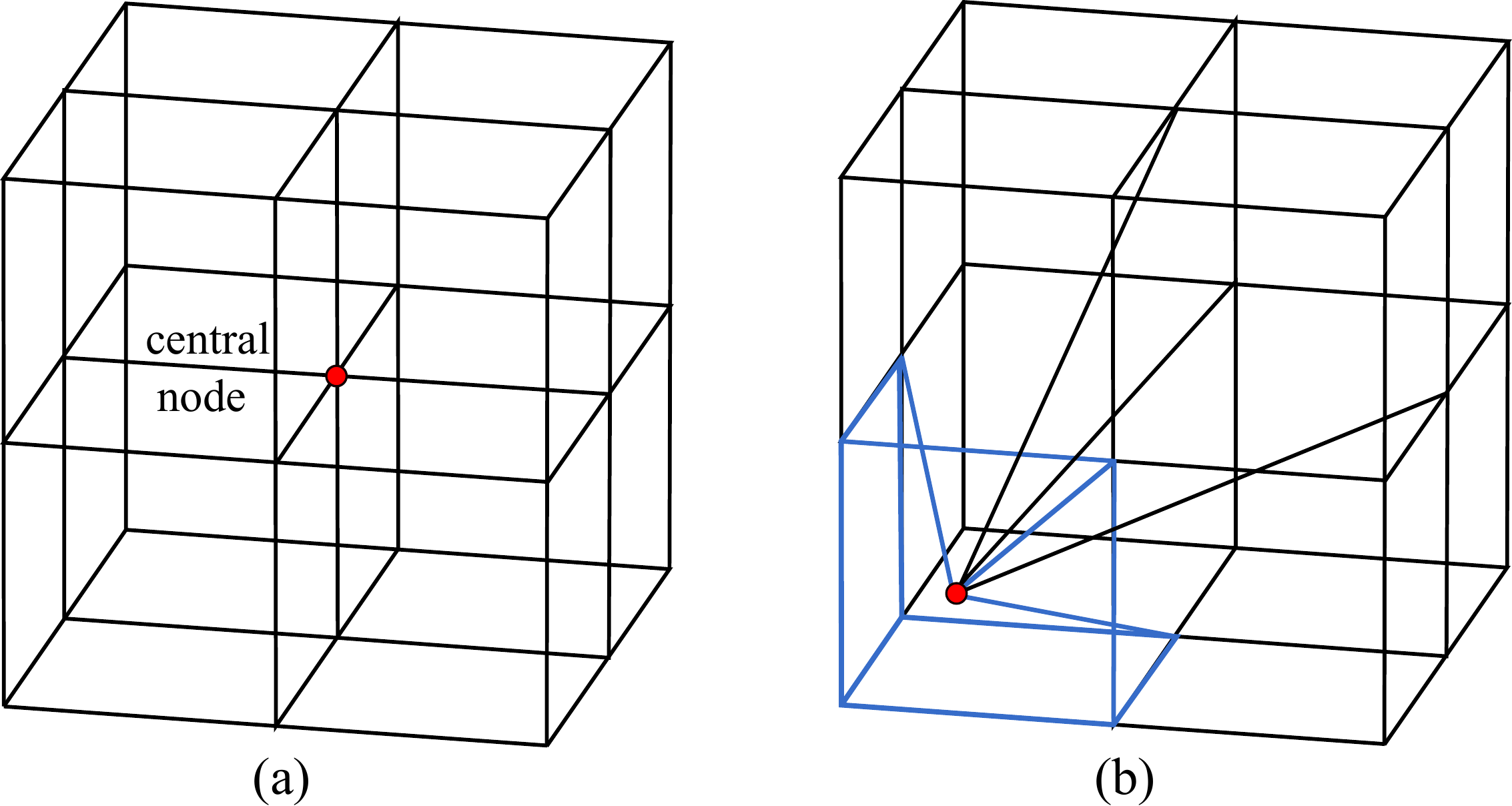}
	\end{center}
	
	\caption{Eight-element (a) regular grid (b) tangled mesh with a concave hex element.}
	\label{Fig_3D_8elem_Mesh}
\end{figure}

\begin{figure}

	\begin{center}
		\includegraphics[width=0.38\textwidth]{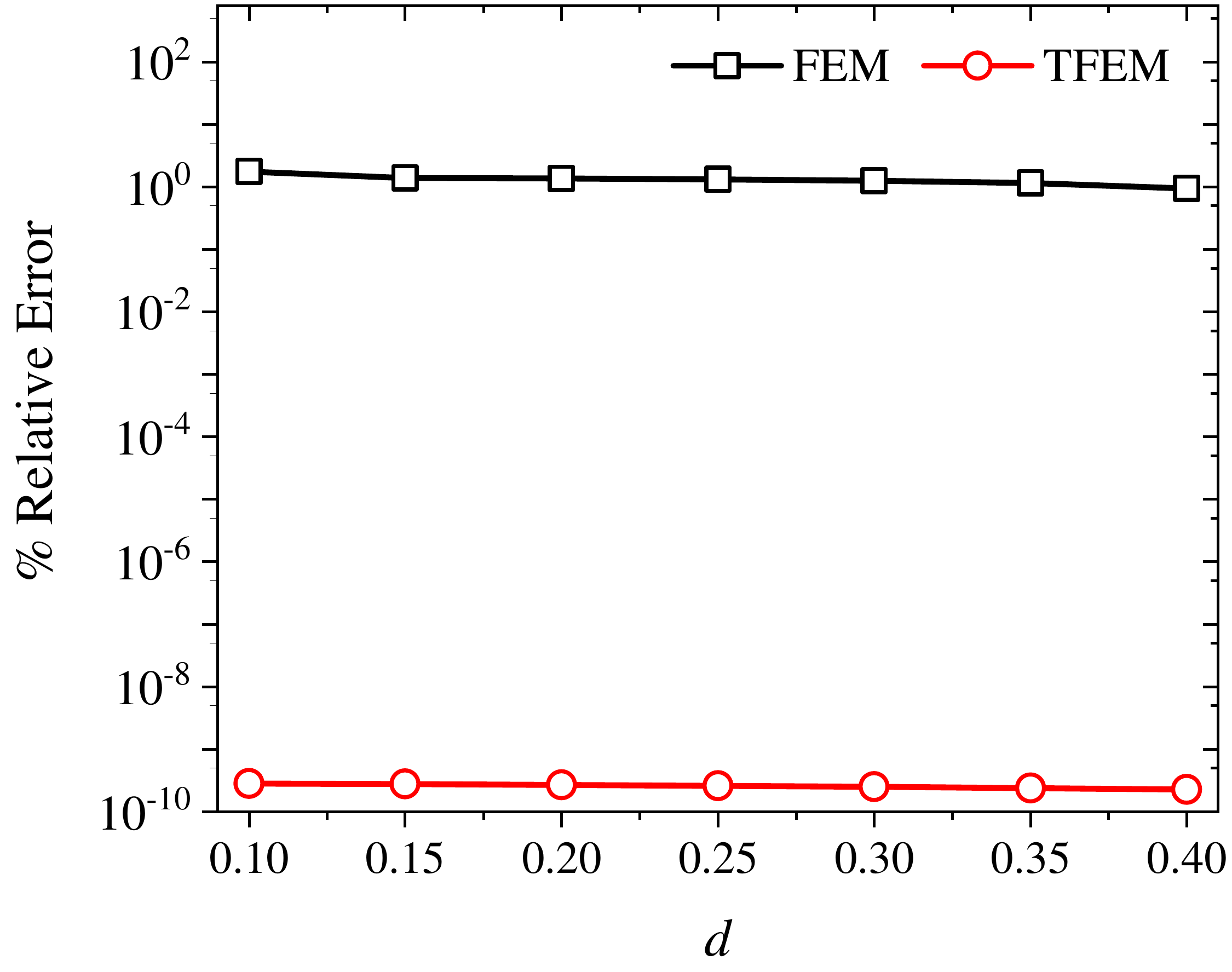}
	\end{center}
	
	\caption{TFEM and FEM errors for eight-element mesh.}
	\label{Fig_3D_8elem_accuracy}
\end{figure}

	Next, consider the a regular uniform grid with 8 hexahedral elements as shown in Fig.~\ref{Fig_3D_8elem_Mesh}a with the central node located at (0.5, 0.5, 0.5). In order to demonstrate TFEM, the mesh can be tangled by moving the central node resulting in one element being concave as shown in Fig.~\ref{Fig_3D_8elem_Mesh}b. The extent of tangling can be varied by defining the position of the central node as $d(0.5, 0.45, 0.4)$ where $ d \in \left[0.1,0.4\right]$. Note that the concave element has only one re-entrant vertex. Fig.~\ref{Fig_3D_8elem_accuracy} compares the relative error in the solution obtained from TFEM and FEM; TFEM provide more accurate results (around $10^{10}$ more accuracy) compared to FEM  while making a few changes in the existing FEM framework.
	
\section{Conclusion} 
\label{Sec4}
This paper demonstrates that the tangled finite element method (TFEM) can provide an effective way to handle concave hexahedral elements in the mesh. TFEM reduces to standard FEM for untangled meshes. To handle tangled elements, the the standard stiffness matrix is modified and an equality condition is added. Numerical experiments illustrate the correctness of the proposed methodology; this is in contrast to standard FEA. 

TFEM solution for various configurations of inverted hex elements and ways to further reduce TFEM error in computationally feasible manner are being investigated. TFEM also extends to nonlinear problems and isogeometric analysis and will be investigated in future. In addition, as tangled meshes can now be handled in TFEM, it would be appropriate to reevaluate mesh generators and optimizers. 

\bibliography{demoref}

\end{document}